\documentclass[12pt,twoside,doublespace]{article}

\usepackage{latexsym,amssymb,amsfonts}
\usepackage{graphicx}

\def\smskip{\par\vskip 5 pt}
\def\QED{\hfill $\Box$\smskip}
\newtheorem{theorem}{Theorem}
\newtheorem{lemma}{Lemma}
\newtheorem{proposition}{Proposition}

\begin{document}

\begin{center}

\vspace{35pt}

{\Large \bf Decentralized Multi-Agent Optimization}

\vspace{5pt}

{\Large \bf  Based on a Penalty Method}

\vspace{5pt}

\vspace{35pt}

{\sc Igor V.~Konnov\footnote{\normalsize E-mail: konn-igor@ya.ru}}

\vspace{35pt}

{\em  Department of System Analysis
and Information Technologies, \\ Kazan Federal University, ul.
Kremlevskaya, 18, Kazan 420008, Russia.}

\end{center}

\vspace{25pt}

\hrule

\vspace{10pt}

{\bf Abstract:} We propose a decentralized penalty method for general convex
constrained multi-agent optimization problems. Each auxiliary penalized
problem is solved approximately with a special parallel descent
splitting method. The method can be implemented in a
computational network where each agent sends information
only to the nearest neighbours. Convergence of the method
is established under rather weak assumptions.
We also describe a specialization of the proposed approach
to the feasibility problem.


\vspace{10pt}

\hrule

\vspace{10pt}

{\em MS Classification:} {65K05, 90C06, 90C25, 68M14, 68W15, 93A14}


\section{Introduction}\label{sc:1}

The custom way of solution of a decision making problem associated with
some complex system consists of collection of all the necessary problem data
in one center and then applying a suitable computational method. For instance,
it can be formulated as an optimization problem and consists
in finding the minimal value of some goal (dis-utility) function $\tilde f$ on
a feasible set $\tilde X$. For brevity, we write
this problem as
\begin{equation} \label{eq:1.1}
 \min \limits _{v \in \tilde X} \to \tilde f (v).
\end{equation}
This means that all the information about the function $\tilde f$ and
set $\tilde X$, which is sufficient for providing efficient computations,
is stored in the central unit, moreover, its computational capacity enables one to
obtain a solution point (or its approximation) within an indicated time period.
However, this situation is not typical for many recent applications related to
large complex systems involving many elements with their private information about
the whole problem and local computational resources. Moreover, the transmission of this
information to the central unit and back is then not suitable since this usually leads to
increasing the data noise and mistakes and to very slow procedures due to various
transmission data delays. In addition, the central unit capacity is smaller essentially than
the total information volume obtained from the whole system. For these reasons,
various decentralized multi-agent procedures become the main direction for solution of
these problems in distributed systems; see e.g. \cite{KPK09,LOF11,PYY13,SFSPP14} and the references therein.

Usually, these problems are then formulated as
optimization problem (\ref{eq:1.1}) where
\begin{equation} \label{eq:1.2}
 \tilde X=\bigcap \limits^{m} _{i=1} X_{i} \ \mbox{and} \ \tilde f (v)=\sum \limits^{m} _{i=1} f_{i}(v),
\end{equation}
$m$ is the number of agents (units) in the system. That is,
the information about the function $f_{i}$ and set $ X_{i}$ is known to the $i$-th agent
and may be unknown even to its neighbours. Besides, it is usually supposed that the system
is connected, i.e. the agents are joined by some transmission links for possible information exchange
so that the system is a connected network. If the problem data are distributed within the complex system
without any preliminary centralized verification
it is natural to suppose that the feasible set $\tilde X$ may appear empty,
and this fact should be also taken into account when creating decentralized solution methods.

The early decomposition methods for large scale optimization problem were mostly oriented on
a significant reduction of the information flows from the central unit to the other units and back,
but the necessity of certain coordination of the whole processes forced one to keep
some central unit; see e.g. \cite{Las70,BLT74,BT89}. Hence, these decomposition methods are not
fully decentralized ones. The modern decentralized optimization methods can be divided into
two main classes. The first class consists of the so-called incremental methods
applied directly to problems of form (\ref{eq:1.1})--(\ref{eq:1.2}); see e.g. \cite{LOF11,DAW12,NO15}
and the references therein.
The second class consists of various primal-dual decomposition methods
applied to their saddle point re-formulation; see e.g. \cite{BPCPE11,LLZ20}
and the references therein.

However, in this paper we intend to develop decentralized penalty methods for
problem (\ref{eq:1.1})--(\ref{eq:1.2}). We recall that the simplest and most popular method for handling
various constraints is the method of smooth penalty functions; see,
e.g. \cite{FM72,Kar75,GK81}, where the original problem is replaced by a
sequence of auxiliary problems with simple constraints.
To the best of our knowledge, this method was not used in multi-agent optimization
since it is not free from flaws.
Firstly, it does not allow one to
find a solution of the original problem with high precision since this
requires very large values of the penalty parameter, but then
finding a solution of penalized problems becomes quite
difficult. Secondly, penalized problems involve binding expressions
for different variables even if the original
problem  is completely decomposable.
Nevertheless, we think these drawbacks do not prevent in fact from application of
penalty-based methods to multi-agent optimization.
We recall that penalty methods are convergent under very general conditions on the problems
in comparison with the other methods, in particular, they are convergent
even if the feasible set is empty. Besides, they are rather stable with respect to
various perturbations and can be applied to non-stationary (limit) problems; see,
e.g. \cite{Kon14b}. Due to essential features of distributed optimization problems
such as utilization of mostly inexact and noisy transmitted information, attaining a
high precision for solutions seems non-realistic for any iterative method, hence penalty methods are in fact
suitable for these problems. Next, several decomposition penalty-based methods were proposed
for large scale optimization problems together with the other decomposition approaches; see, e.g.
\cite{Raz67,Mau71,Raz75,Umn75}. These methods also involved some central unit due to the necessity
of certain coordination of the whole iterative process. Rather recently, some other
decomposition penalty method was proposed in \cite{Kon19c}, which is based on a special approximation of
auxiliary penalized problems and a descent splitting method. In this paper we combine this technique
and peculiarities of the penalized formulation of the multi-agent optimization problem, which leads to a
decentralized multi-agent penalty process. We prove convergence of the proposed penalty method
under rather weak assumptions and describe its information exchange scheme.
Besides, we describe a specialization of the proposed approach
to the feasibility problem and give results of preliminary computational experiments,
which showed rather satisfactory and stable convergence.

We outline now briefly the further organization of the paper.
In Section \ref{sc:2}, we recall some auxiliary properties and facts from
the theory of convex optimization.
 In Section \ref{sc:3}, we  give a re-formulation of the original problem
(\ref{eq:1.1})--(\ref{eq:1.2}) and substantiate a general penalty method.
In Section \ref{sc:4}, we  present a decentralized two-level penalty method
with approximate solution of each penalized problem and prove its convergence.
Implementation issues of the method are discussed in
Section \ref{sc:5}. In Section \ref{sc:6}, we  describe an application of the proposed approach
to the feasibility problem.
Section \ref{sc:7} gives examples of preliminary calculations
of the proposed method and some other basic decomposition and multi-agent
methods on test problems. Section \ref{sc:8} contains some conclusions.


\section{ Auxiliary properties }\label{sc:2}

This section presents some results from the theory of
convex optimization that will be used in the next sections. Let us consider first the
optimization problem
\begin{equation} \label{eq:2.1}
 \min \limits _{x \in X} \to \mu (x),
\end{equation}
for some  function $\mu$ and set $X$,
the set of its solutions is denoted by $X^{*}(\mu)$, and the optimal
function value by $\mu^{*}$, i.e.
$$
\mu^{*} = \inf \limits_{ x \in X} \mu(x).
$$
If the set $X^{*}(\mu)$ is bounded, then the function $\mu$ has a linear
minorant on the set $X$; see \cite[Ch. IV, \S 2, Theorem
18]{Vas11}. Here we present this result in a somewhat modified
format.


\begin{proposition} \label{pro:2.1} Let $X$ be a
convex and closed set in $ \mathbb{R}^{N}$ and $\mu: X \to
\mathbb{R}$ a continuous convex function. If the set
$X^{*}(\mu)$ is non-empty and bounded, then for any point $x^{*} \in
X^{*}(\mu)$ there are a bounded set $U\supset X^{*}(\mu)$ and
a number $\sigma>0$ such that
$$
 \mu(x)-\mu^{*}  \geq \sigma \|x-x^{*} \|, \quad \forall x\in X
\setminus U.
$$
\end{proposition}

We give an optimality condition in the additive case when
\begin{equation} \label{eq:2.2}
\mu(x)=\mu_{1}(x)+\mu_{2}(x),
\end{equation}
where $\mu_{2}: \mathbb{R}^{N} \to \mathbb{R}$ is a smooth function; see
\cite[Proposition 2.2.2]{Pan85} and \cite[Proposition 1]{KS17}.
We will use the following basic assumptions.
 \begin{enumerate}
 \item[{\bf (A1)}] $ X $ is a nonempty, convex, and closed set in $ \mathbb{R}^{N} $.
 \item[{\bf (A2)}] $\mu_{1} : \mathbb{R}^{N} \to \mathbb{R}$ is a convex function, $\mu_{2} :
\mathbb{R}^{N} \to \mathbb{R}$ is a  smooth convex function.
 \end{enumerate}


\begin{proposition} \label{pro:2.2}  Let conditions (A1)--(A2) be satisfied.  Then problem
(\ref{eq:2.1})--(\ref{eq:2.2}) is equivalent to the mixed
variational inequality (MVI for short): Find a point $x^{*} \in X$ such that
\begin{equation} \label{eq:2.3}
   \left[\mu_{1}(x)-\mu_{1}(x^{*}) \right]
   + \langle  \mu'_{2}(x^{*}), x-x^{*} \rangle   \geq 0 \quad \forall x \in X.
\end{equation}
\end{proposition}

Let's fix a number $\alpha >0$ and consider the auxiliary  optimization problem:
\begin{equation} \label{eq:2.4}
\begin{array}{c}
\displaystyle \min \limits_{z \in X} \rightarrow \left\{
{\mu_{1}(z)+\langle \mu'_{2}(x), z \rangle + (2\alpha )^{-1}
\|z-x\|^{2}} \right\}
\end{array}
\end{equation}
for some point $x \in X$. Under the assumptions made, the goal
function in (\ref{eq:2.4}) is continuous and strongly convex, so
 problem (\ref{eq:2.4}) has a unique solution which we denote by
$y_{\alpha}(x)$, thus defining a single-valued mapping
$x\mapsto y_{\alpha}(x)$. Instead of problem (\ref{eq:2.4}) it will be convenient
to use also its equivalent formulation in the form of a MVI.

\begin{lemma} \label{lm:2.1} Let conditions (A1)--(A2) be satisfied.
The point $y_{\alpha}(x) \in X$ is a solution to problem (\ref{eq:2.4})
if and only if it satisfies the condition
\begin{equation} \label{eq:2.5}
\begin{array}{c}
\displaystyle \mu_{1}(z)-\mu_{1}(y_{\alpha}(x))+\alpha ^{-1} \langle
y_{\alpha}(x)-x, z-y_{\alpha}(x)\rangle  \\
\displaystyle + \langle \mu'_{2}(x), z-y_{\alpha}(x)\rangle
\geq 0 \quad \forall z \in X.
\end{array}
\end{equation}
\end{lemma}
{\bf Proof.} Problem (\ref{eq:2.4}) is written equivalently
 in the form:
$$
\min \limits_{y \in X} \rightarrow f_{1}(y)+f_{2}(y),
$$
where $f_{1}(z)=\mu_{1}(z)+\langle \mu'_{2}(x), z \rangle$ and $f_{2}(z)=(2\alpha )^{-1} \| z-x\|^{2}$.
Note that the functions $f_{1}$ and $f_{2}$ are
convex, and the function $f_{2}$ is differentiable.
According to Proposition \ref{pro:2.2}, this problem is equivalent to the MVI:
$$
f_{1}(z)-f_{1}(y_{\alpha}(x))
+ \langle f'_{2}(y_{\alpha}(x)), z-y_{\alpha}(x)\rangle
 \geq 0 \quad \forall z \in X,
$$
which obviously coincides with (\ref{eq:2.5}). \QED

Now we get a few basic  properties of the mapping  $x\mapsto
y_{\alpha}(x)$.


\begin{proposition} \label{pro:2.3}  Let conditions (A1)--(A2) be satisfied.  Then
the following statements are true.

(a) The set of fixed points of the mapping $x\mapsto y_{\alpha}(x)$
 coincides with the set of solutions of problem (\ref{eq:2.1})--(\ref{eq:2.2});

(b)  The mapping $x\mapsto y_{\alpha}(x)$ is continuous on $X$;

(c) For any point $z \in X$ it holds that
\begin{equation} \label{eq:2.6}
\begin{array}{l}
\mu_{1}(y_{\alpha}(x))-\mu_{1}(z) + \langle \mu'_{2}(y_{\alpha}(x)), y_{\alpha}(x)-z\rangle\\
\displaystyle \leq \langle \mu'_{2}(y_{\alpha}(x))-\mu'_{2}(x), y_{\alpha}(x)-z\rangle
+\alpha ^{-1} \langle y_{\alpha}(x)-x, z-y_{\alpha}(x)\rangle .
\end{array}
\end{equation}
\end{proposition}
{\bf Proof.}
If $ x^{*} = y_{\alpha}(x^{*} )$, then
(\ref{eq:2.5}) implies $x^{*} \in X^{*}(\mu)$. Conversely, let
 $x^{*}$ solve MVI (\ref{eq:2.3}), but $ x\neq y_{\alpha}(x)$.
 Then setting $z=x$ in (\ref{eq:2.5}) gives
$$
\mu_{1}(x)-\mu_{1}(y_{\alpha}(x))+\langle \mu'_{2}(x), x-y_{\alpha}(x)\rangle
\geq  \alpha
^{-1}\|y_{\alpha}(x)-x\|^{2}>0,
$$
which is a contradiction. Part (a) is true. To prove (b), take
arbitrary $x'$, $x'' \in X$ and set $y'=y_{\alpha}(x')$
and $y''=y_{\alpha}(x'')$ for brevity. Then from (\ref{eq:2.5}) it follows
that
$$
\mu_{1}(y'')-\mu_{1}(y') +\langle \mu'_{2}(x')+\alpha^{-1}(y'-x'), y''-y'\rangle \geq 0
$$
and
$$
\mu_{1}(y')-\mu_{1}(y'') +\langle \mu'_{2}(x'')+\alpha^{-1}(y''-x''), y'-y''\rangle \geq 0.
$$
Summing these inequalities gives
$$
\langle \mu'_{2}(x')-\mu'_{2}(x'')-\alpha^{-1}(x'-x''),y''-y' \rangle
 \geq \alpha ^{-1}\|y''-y'\|^{2},
$$
hence
$$
\|\mu'_{2}(x')-\mu'_{2}(x'')\| +\alpha
^{-1}\|x'-x''\| \geq \alpha
^{-1}\|y''-y'\|.
$$
This means the mapping $x\mapsto y_{\alpha}(x)$ is
continuous and part (b) is true. To prove (c), we again use
(\ref{eq:2.5}) and obtain
\begin{eqnarray*}
\displaystyle && \mu_{1}(y_{\alpha}(x))-\mu_{1}(z) + \langle \mu'_{2}(y_{\alpha}(x)), y_{\alpha}(x)-z\rangle \\
&& = \mu_{1}(y_{\alpha}(x))-\mu_{1}(z) + \langle \mu'_{2}(x)+\alpha ^{-1} (y_{\alpha}(x)-x), y_{\alpha}(x)-z\rangle  \\
&& + \langle \mu'_{2}(y_{\alpha}(x))-\mu'_{2}(x)-\alpha ^{-1} (y_{\alpha}(x)-x), y_{\alpha}(x)-z\rangle \\
&& \leq \langle \mu'_{2}(y_{\alpha}(x))-\mu'_{2}(x)-\alpha ^{-1} (y_{\alpha}(x)-x), y_{\alpha}(x)-z\rangle,
\end{eqnarray*}
which gives (\ref{eq:2.6}).
\QED

Next, we recall that the  iterate
\begin{equation} \label{eq:2.7}
x^{k+1}=y_{\alpha}( x^{k}), \ k=0,1,2,\ldots,
\end{equation}
corresponds to the well-known forward-backward
splitting method; see \cite{LM79,Gab83}.
We intend to utilize its descent properties. We need the basic inequality for
a function having the Lipschitz continuous gradient; see \cite[Ch. III, Lemma 1.2]{DR68}.


\begin{proposition} \label{pro:2.4} Let  the gradient of a
function $\varphi : \mathbb{R}^{N} \to \mathbb{R}$  satisfies the Lipschitz
condition with constant $L_{\varphi}$ on a convex set $X$. Then
$$
\varphi(x'') \leq \varphi(x')+\langle \varphi'(x'),x''-x' \rangle +0.5L_{\varphi} \|x''-x'\|^{2}  \quad \forall x',x''\in X.
$$
\end{proposition}

We now somewhat modify the assumptions in (A2).
 \begin{enumerate}
 \item[{\bf (A2$'$)}] $\mu_{1} : \mathbb{R}^{N} \to \mathbb{R}$ is a convex function, $\mu_{2} :
\mathbb{R}^{N} \to \mathbb{R}$ is a  smooth convex function, its gradient satisfies the Lipschitz
condition with constant $L_{\mu_{2}}$ on the set $X$.
 \end{enumerate}

\begin{lemma} \label{lm:2.2} Let conditions (A1) and (A2$\,'$) be satisfied.
If
\begin{equation} \label{eq:2.9}
\alpha \leq 1/(\beta+0.5L_{\mu_{2}})
\end{equation}
for some $\beta \in (0,1)$, then iterate (\ref{eq:2.7}) yields
\begin{equation} \label{eq:2.10}
\mu(x^{k+1}) \leq \mu(x^{k})-\beta \|x^{k+1}- x^{k} \| ^{2}.
\end{equation}
\end{lemma}
{\bf Proof.} By definition,
$$
\mu(x^{k+1})-\mu(x^{k})= \mu_{1}(x^{k+1})-\mu_{1}(x^{k}) +\mu_{2}(x^{k+1})-\mu_{2}(x^{k}).
$$
Applying Proposition \ref{pro:2.4} to $\mu_{2}$, we have
$$
\mu_{2}(x^{k+1})-\mu_{2}(x^{k}) \leq \langle \mu'_{2}(x^{k}),x^{k+1}- x^{k} \rangle +0.5L_{\mu_{2}} \|x^{k+1}- x^{k}\|^{2}.
$$
Combining these relations with (\ref{eq:2.5}) where $x=z=x^{k}$, we obtain
\begin{eqnarray*}
\displaystyle && \mu(x^{k+1})-\mu(x^{k}) \leq \mu_{1}(x^{k+1})-\mu_{1}(x^{k})
+\langle \mu'_{2}(x^{k}),x^{k+1}- x^{k} \rangle +0.5L_{\mu_{2}} \|x^{k+1}- x^{k}\|^{2}  \\
&& \leq  -\alpha ^{-1} \|x^{k+1}- x^{k} \| ^{2}+ 0.5L_{\mu_{2}} \|x^{k+1}- x^{k}\|^{2} \\
&& = -(\alpha ^{-1}-0.5L_{\mu_{2}}) \|x^{k+1}- x^{k}\|^{2}.
\end{eqnarray*}
Due to (\ref{eq:2.9}), this inequality gives (\ref{eq:2.10}).
\QED

We also recall that a function $\varphi : \mathbb{R}^{N} \to \mathbb{R}$  is coercive on a set $X$
if
$$
\varphi (x) \to +\infty \quad \mbox{as} \quad \|x\| \to \infty, \ x \in X.
$$
After adding the coercivity assumption for the function $\mu$ we immediately obtain the basic
convergence properties for the splitting method (\ref{eq:2.7}) from Lemma \ref{lm:2.2}.


\begin{proposition} \label{pro:2.5}  Let conditions (A1) and (A2$\,'$) be satisfied,
 the function $\mu$ be coercive on the set $X$, and let the
sequence $\{x^{k}\}$ be generated in accordance with rules (\ref{eq:2.7}) and (\ref{eq:2.9})
for some $\beta \in (0,1)$.  Then
the  sequence $\{x^{k}\}$  has limit points, all these limit points are solutions of
problem (\ref{eq:2.1})--(\ref{eq:2.2}), besides,
$$
\lim \limits_{k\rightarrow \infty }\|x^{k+1}- x^{k}\|=0
$$
and
$$
\lim \limits_{k\rightarrow \infty }\mu (x^{k})=\mu^{*}.
$$
\end{proposition}


\section{The basic problem re-formulation and a penalty method}\label{sc:3}

We now present a re-formulation of the original problem
(\ref{eq:1.1})--(\ref{eq:1.2}). We write
this problem as
\begin{equation} \label{eq:3.1}
 \min \limits _{x \in D} \to f(x) = \sum \limits^{m} _{i=1} f_{i}(x_{i}),
\end{equation}
where $x=(x_{i})_{i=1, \ldots, m} \in \mathbb{R}^{N}$, i.e. $x^{\top}=(x^{\top}_{1}, \dots,x^{\top}_{m})$,
$x_{i}= (x_{i1}, \dots,x_{in})^{\top}$ for $i=1, \dots, m$, $N=mn$,
\begin{equation} \label{eq:3.2}
 D = X \bigcap Y, \ X=X_{1}\times \dots \times X_{m}=\prod \limits_{i=1}^{m} X_{i},
\end{equation}
the set $Y$ describes the topology of the communication network. For instance, the set
$$
Y'=\left\{x \in \mathbb{R}^{N} \ | \ x_{i}=x_{i+1}, \ i=1, \dots, m-1  \right\}
$$
gives the minimal connected graph topology (chain), whereas
$$
Y''=\left\{x \in \mathbb{R}^{N} \ | \ x_{i}=x_{j}, \ i,j=1, \dots, m, \ i \neq j  \right\},
$$
gives the maximal (full) graph. In principle, we can take any suitable variant between $Y'$ and  $Y''$.
In order to both increase the communication reliability and reduce the transmission flows we take the
following set
\begin{equation} \label{eq:3.3}
Y=\left\{x \in \mathbb{R}^{N} \ | \ x_{i}=x_{i+1}, \ i=1, \dots, m-1, x_{m}=x_{1}  \right\},
\end{equation}
which corresponds to the simplest cycle in the system.
Here each unit receives information only from two neighbours.
Clearly, the unit numbering can be chosen arbitrary.
However, in the case where the set $\tilde X$ (hence $D$) may be empty,
we should change the formulation.
The constraints of the set $Y$ will be taken into account by the penalty function
\begin{equation} \label{eq:3.4}
 p (x) = (2\tau )^{-1} \| Ax\|^{2}=(2\tau )^{-1}\sum \limits^{m} _{i=1} \| A_{i}x\|^{2},
\end{equation}
where
$$
A= \left( {
\begin{array}{cccccc}
I & -I &   \Theta & \dots & \Theta & \Theta \\
\Theta & I & -I & \dots & \Theta & \Theta\\
\dots & \dots &  \dots & \dots & \dots & \dots \\
-I & \Theta & \Theta &  \dots & \Theta & I
\end{array}
} \right)=
\left( {
\begin{array}{c}
A_{1} \\
A_{2} \\
\dots \\
A_{m}
\end{array}
} \right),
$$
$I$ is the $n \times n$ unit matrix, $\Theta$ is the $n \times n$ zero matrix,
$A_{i}$ is the corresponding $n \times nm$ sub-matrix of $A$ for $i=1, \dots, m$, and
$\tau > 0$ is a fixed scaling parameter.
Let
$$
p^{*} = \inf \limits_{ x \in X} p(x).
$$
Then we can take the more general (sequential) optimization problem:
\begin{equation} \label{eq:3.4a}
 \min \limits _{x \in  X^{*}(p)} \to f(x),
\end{equation}
and denote by $\tilde D$ its solution set.
Clearly, problem (\ref{eq:3.4a}) coincides with (\ref{eq:3.1})--(\ref{eq:3.3})
if $D \neq \varnothing$ and $p^{*} = 0$.

In what follows, we will use the following basic assumptions.
\begin{enumerate}
\item[{\bf (B1)}] The set $ X^{*}(p)$ is nonempty, $X_{i}$ is a  convex and closed set in
$ \mathbb{R}^{n}$ for $i=1, \dots, m$.
\item[{\bf (B2)}] $f : \mathbb{R}^{N} \to \mathbb{R}$ is a coercive function on a set $X$,
$f_{i} : \mathbb{R}^{n} \to \mathbb{R}$ is a convex
function for $i=1, \dots, m$.
\end{enumerate}

Under the above assumptions
problem (\ref{eq:3.4a})  has a solution, i.e. the sets $X^{*}(f)$ and $\tilde D$ are nonempty
and  bounded. Set
$$
f^{**} = \inf \limits_{ x \in X^{*}(p)} f(x), \  f^{*} = \inf \limits_{ x \in D} f(x), \ \mbox{and} \ \phi^{*} = \inf \limits_{ x \in X} f(x).
$$
It is clear that
$$
\phi^{*} = \sum \limits^{m} _{i=1} \phi^{*}_{i} \ \mbox{where} \
 \phi^{*}_{i} = \inf \limits_{ x_{i} \in X_{i}} f_{i}(x_{i}), \ i=1, \dots, m.
$$

For a fixed vector $e= (\varepsilon_{1}, \dots,\varepsilon_{m})^{\top}$ of positive penalty
parameters we can define the auxiliary function
\begin{equation} \label{eq:3.5}
 \varphi(x,e)=h(x,e)+p(x), \ h(x,e)=\langle  e, F(x) \rangle=\sum \limits^{m} _{i=1} \varepsilon_{i} f_{i}(x_{i}),
\end{equation}
where $F(x)= (f_{1}(x_{1}), \dots,f_{m}(x_{m}))^{\top}$.
So, the initial problem (\ref{eq:3.4a})
is replaced by a sequence of auxiliary problems of the form
\begin{equation} \label{eq:3.6}
 \min \limits _{x \in X} \to \varphi(x,e).
\end{equation}
with separate constraints. The custom penalty approach utilizes one scalar
penalty parameter (see \cite{Kar75}), but here each agent can manage his/her own
penalty parameter, which seems more natural for decomposable systems.

We denote by $z(e)$ any solution of problem (\ref{eq:3.5})--(\ref{eq:3.6}).
Our first goal is to prove that the trajectory $\{z(e)\}$ tends in some
sense to  a solution of problem (\ref{eq:3.1})--(\ref{eq:3.3}) as $e
\to \mathbf{0}$. Also, for brevity, we set $z^{s}=x(e^{s})$
for any sequence $\{e^{s}\}$.


\begin{theorem} \label{thm:3.1}
Suppose that assumptions (B1)--(B2) are fulfilled, the sequence
$\{e^{s}\}$ satisfies the conditions:
\begin{equation} \label{eq:3.7}
\{\varepsilon^{s}_{i}\} \searrow 0, \ i=1, \dots, m, \
\lim \limits_{s \rightarrow \infty } (\varepsilon^{s}_{i}/\varepsilon^{s}_{j})=1,
\ \forall i \neq j.
\end{equation}
Then:

(i) Problem (\ref{eq:3.6}) has a solution for each positive vector $e$;

(ii) Each sequence $\{z^{s}\}$ of solutions of (\ref{eq:3.6}) has
limit points and all these limit points are solutions of problem (\ref{eq:3.4a}).
 \end{theorem}
{\bf Proof.} Due to (B2) each function $f_{i}$ is  coercive  on a set $X_{i}$, hence
the function $\varphi$ is  coercive  on $X$ and problem (\ref{eq:3.6}) has a solution, i.e.
part (i) is true. Therefore, the above penalty method is well-defined.

Next, take any point $x^{*} \in \tilde D$. Then from the definition we have
\begin{equation} \label{eq:3.8}
h(z^{s},e^{s}) \leq \varphi (z^{s},e^{s}) \leq \varphi (x^{*} ,e^{s})=h(x^{*},e^{s})+p(x^{*})
  \leq h(x^{*},e^{s})+p(z^{s}).
\end{equation}
Assume that $ \{ \|z^{s}\| \} \to \infty$. Then there exists at least one index $j$ such that
$ \{ f_{j}(z^{s}_{j}) \} \to +\infty$. From (\ref{eq:3.8}) it follows that
$$
h(z^{s},e^{s})-h(x^{*},e^{s}) =\langle  e^{s}, F(z^{s}) -F(x^{*}) \rangle \leq 0
$$
Dividing both the sides  on $\varepsilon^{s}_{j}$ we obtain
$$
f_{j}(z^{s}_{j}) \leq f_{j}(x^{*}_{j})
+\sum \limits _{i \neq j} (\varepsilon^{s}_{i}/\varepsilon^{s}_{j})( f_{i}(x^{*}_{i})-\phi^{*}_{i}),
$$
taking now the limit $s \rightarrow \infty $ gives a contradiction. Hence, the
sequence $\{z^{s}\}$ is bounded and has
limit points. Let $\bar z$ be an arbitrary limit point for $\{z^{s}\}$,
i.e.
$$
\bar z = \lim \limits_{{l} \rightarrow \infty } z^{s_{l}}.
$$
Then clearly $\bar z \in X$. From (\ref{eq:3.8}) it follows that
$$
0 \leq p(z^{s}) \leq  p^{*}+\sum \limits^{m} _{i=1} \varepsilon^{s}_{i}( f_{i}(x^{*}_{i})-\phi^{*}_{i}),
$$
taking the limit $s=s_{l} \rightarrow \infty $ gives $p(\bar z) \leq p^{*}$,
hence $p(\bar z)=p^{*}$. This means that all the limit points of $\{z^{s}\}$
belong to $X^{*}(p)$. Again from (\ref{eq:3.8}) we have
$$
\sum \limits^{m} _{i=1} (\varepsilon^{s}_{i}/\varepsilon^{s}_{j})f_{i}(z^{s}_{i}) \leq  \sum \limits^{m} _{i=1} (\varepsilon^{s}_{i}/\varepsilon^{s}_{j}) f_{i}(x^{*}_{i}).
$$
Taking the limit $s=s_{l} \rightarrow \infty $ gives
$$
f(\bar z)=\sum \limits^{m} _{i=1} f_{i}(\bar z_{i}) \leq  \sum \limits^{m} _{i=1} f_{i}(x^{*}_{i})=f^{**}.
$$
Hence, part (ii) is also true. \QED


\section{Two-level penalty method}\label{sc:4}

The goal function $\varphi(x,e)$ of problem (\ref{eq:3.6}) still involves coupled variables
 that prevents from direct application of decentralized control schemes. For this reason,
we intend to find an approximate solution of each penalized problem by using the forward-backward
splitting method (\ref{eq:2.7}). More precisely, we fix a number $\alpha >0$
and define the point $y_{\alpha,e}(x)$ as a unique solution of
the optimization problem:
\begin{equation} \label{eq:4.1}
\begin{array}{c}
\displaystyle \min \limits_{z \in X} \rightarrow \left\{
{h(z,e)+\langle p'(x), z \rangle + (2\alpha )^{-1}
\|z-x\|^{2}} \right\}
\end{array}
\end{equation}
for some point $x \in X$; cf. (\ref{eq:2.4}). Then iterate (\ref{eq:2.7}) is re-written as follows:
\begin{equation} \label{eq:4.2}
x^{k+1}=y_{\alpha,e}( x^{k}), \ k=0,1,2,\ldots,
\end{equation}
its convergence will follow from Lemma \ref{lm:2.2} and
Proposition \ref{pro:2.5}, but we have to evaluate
the Lipschitz constant $L_{p}$ of the gradient $p'(x)$.

\begin{lemma} \label{lm:4.1}
It holds that $L_{p} \leq 4/\tau$.
\end{lemma}
{\bf Proof.} It follows from (\ref{eq:3.4}) that
$$
p (x) = (2\tau )^{-1} \| Ax\|^{2}=(2\tau )^{-1} \langle Sx, x \rangle,
$$
where
$$
S=A^{\top}A= \left( {
\begin{array}{cccccc}
2I & -I &   \Theta & \dots & \Theta & -I \\
-I & 2I & -I & \dots & \Theta & \Theta\\
\Theta & -I & 2I & \dots & \Theta & \Theta\\
\dots & \dots &  \dots & \dots & \dots & \dots \\
\Theta & \Theta & \Theta &  \dots & 2I & -I\\
-I & \Theta & \Theta &  \dots & -I & 2I
\end{array}
} \right),
$$
hence $L_{p}=\| S\|/\tau$. Set
$$
S'=  \left( {
\begin{array}{rrrrrr}
2 & -1 &   0 & \dots & 0 & -1 \\
-1 & 2 & -1 & \dots & 0 & 0\\
0 & -1 & 2 & \dots & 0 & 0\\
\dots & \dots &  \dots & \dots & \dots & \dots \\
0 & 0 & 0 &  \dots & 2 & -1\\
-1 & 0 & 0 &  \dots & -1 & 2
\end{array}
} \right),
$$
then
$$
S=  S' \otimes I,
$$
where $\otimes$ denotes the Kronecker product of matrices,
hence the eigenvalues of $S$ and $S'$ coincide; see Theorem 3 in \cite[Ch.XII]{Bel97}.
From the Gershgorin theorem (see Theorem 5 in \cite[Ch.XIV]{Gan66})
we obtain that the maximal eigenvalue of $S'$ is not greater than four
and the result follows.
\QED

We will describe the two-level decomposable penalty method, which uses
approximate solutions of problems (\ref{eq:3.6}).

\medskip
\noindent {\bf Method (DPM).} Choose a point  $u^{0} \in X$,
a sequence of positive vectors $\{e^{s}\}$
and a sequence of positive numbers $\{\theta_{s}\}$. Fix numbers $\beta \in (0,1)$
and $\tau  \geq 1$
and choose a number $\alpha \in (0,\alpha']$ where $\alpha' = 1/(\beta+2/\tau)$.

   At the $s$-th stage, $s=1,2,\ldots$, we have a point
$u^{s-1}\in X$ and parameters $e^{s}$ and $\theta_{s}$. Applying
iterate  (\ref{eq:4.2}) with the starting point $x^{0}=u^{s - 1}$ and $e=e^{s}$,
we obtain the point $x^{k+1}$ such that
\begin{equation} \label{eq:4.4}
\|x^{k+1}- x^{k} \| \leq \theta_{s},
\end{equation}
and set $u^{s}= x^{k+1}$.
\medskip

We now obtain the basic convergence statement for (DPM).


\begin{theorem} \label{thm:4.1}
Suppose that assumptions (B1)--(B2) are fulfilled, the sequence
$\{e^{s}\}$ satisfies the conditions in (\ref{eq:3.7}), besides,
\begin{equation} \label{eq:4.5}
\lim \limits_{s \rightarrow \infty }(\theta_{s}/\varepsilon^{s}_{j}) = 0, \ \forall j=1, \dots, m.
\end{equation}
Then:

(i) The number of iterations at each stage of Method (DPM) is finite;

(ii) Each sequence $\{u^{s}\}$ has
limit points and all these limit points are solutions of problem (\ref{eq:3.4a}).
 \end{theorem}
{\bf Proof.} First we note that part (i) is true due to
Proposition \ref{pro:2.5} since the choice of $\alpha$ and $\tau$ in (DPM)
provides (\ref{eq:2.9}) with
respect to any problem (\ref{eq:3.6}) where $\mu_{1}(x)=h(x,e)$ and
$\mu_{2}(x)=p(x)$ because of Lemma \ref{lm:4.1}.

Therefore, the sequence $\{u^{s}\}$ is well defined.
Using (\ref{eq:2.6}) with $\mu_{1}(x)=h(x,e)$ and
$\mu_{2}(x)=p(x)$, Lemma \ref{lm:4.1}, and (\ref{eq:4.4}) we obtain
\begin{equation} \label{eq:4.6}
\langle  e^{s}, F(u^{s}) -F(x^{*}) \rangle \leq (4\tau^{-1}+\alpha ^{-1}) \theta_{s} \| u^{s} -x^{*}\| +
p(x^{*})- p(u^{s})
\end{equation}
for any $x^{*} \in \tilde D$. For
brevity, fix $L=(4\tau^{-1}+\alpha ^{-1})$. It follows that
$$
\tilde \varepsilon_{s} (f(u^{s}) -\phi^{*}) \leq L \theta_{s} (\| u^{s} - \tilde x\| + \| \tilde x -x^{*}\|)
+ \langle  e^{s}, F(\tilde x) -F(x^{*}) \rangle
$$
for any $\tilde x \in X$ such that $f(\tilde x) =\phi^{*}$ where
$$
\tilde \varepsilon_{s}=\min \limits_{i=1, \dots, m} \varepsilon^{s}_{i}.
$$
Assume that $ \{ \|u^{s} \| \} \to +\infty$. Then, taking into account
Proposition \ref{pro:2.1} and dividing the above relation by $\tilde  \varepsilon_{s} \| u^{s} - \tilde x\|$, we obtain
$$
0 < \sigma \leq L (\theta_{s}/\tilde  \varepsilon_{s}) (1 + \| \tilde x -x^{*}\|/\| u^{s} - \tilde x\|)
+ \langle  e^{s}, F(\tilde x) -F(x^{*}) \rangle /(\tilde  \varepsilon_{s} \| u^{s} - \tilde x\|)
$$
for $s$ large enough. It is clear that there exists a subsequence $\{s_{l}\}$
and a fixed index $j$ such that
$$
\tilde \varepsilon_{s_{l}}=\varepsilon^{s_{l}}_{j}.
$$
Taking now the limit $s=s_{l} \rightarrow \infty $ gives
the contradiction  $0<\sigma \leq 0$  due to (\ref{eq:3.7}) and (\ref{eq:4.5}).
Therefore, the sequence $\{u^{s}\}$ is bounded and
has limit points.

Let $\bar u$ be an arbitrary limit point for $\{u^{s}\}$,
i.e.
$$
\bar u = \lim \limits_{{l} \rightarrow \infty } u^{s_{l}}.
$$
Then clearly $\bar u \in X$. From (\ref{eq:4.6}) it follows that
$$
0 \leq p(u^{s}) \leq  p^{*}+\sum \limits^{m} _{i=1} \varepsilon^{s}_{i}( f_{i}(x^{*}_{i})-\phi^{*}_{i})
+ L \theta_{s}\| u^{s} -x^{*}\|.
$$
Taking here the limit $s=s_{l} \rightarrow \infty $ gives $p(\bar u) \leq p^{*}$,
hence $p(\bar u)=p^{*}$. This means that all the limit points of $\{u^{s}\}$
belong to $X^{*}(p)$.

Again from (\ref{eq:4.6}) we have
$$
\sum \limits^{m} _{i=1} (\varepsilon^{s}_{i}/\varepsilon^{s}_{j})f_{i}(u^{s}_{i}) \leq  \sum \limits^{m} _{i=1} (\varepsilon^{s}_{i}/\varepsilon^{s}_{j}) f_{i}(x^{*}_{i})+ L (\theta_{s}/\varepsilon^{s}_{j})\| u^{s} -x^{*}\|.
$$
Taking the limit $s=s_{l} \rightarrow \infty $ and using (\ref{eq:3.7}) and (\ref{eq:4.5}) gives
$$
f(\bar u)=\sum \limits^{m} _{i=1} f_{i}(\bar u_{i}) \leq  \sum \limits^{m} _{i=1} f_{i}(x^{*}_{i})=f^{**}.
$$
Hence, part (ii) is also true. \QED


\section{Implementation issues} \label{sc:5}

In this section we describe a decentralized implementation
of the presented Method (DPM), where each agent (or unit) receives information only from
his/her two closest neighbours and the topology of the communication network is the simplest cycle.
Next, each $i$-th agent keeps his/her private information about the set $X_{i}$ and
function $f_{i}$, which are in general unknown to the others.
Also, the $i$-th agent tells the state $x^{k}_{i}$ to the closest neighbours after its calculation
and may send in principle some other short signals about the calculation process.

 At the beginning of the calculation procedure
 each $i$-th agent has a fixed scaling parameter $\tau  \geq 1$, a proper step-size $\alpha \in (0,\alpha']$,
sequences of positive numbers
 $\{\varepsilon^{s}_{i}\}$ and $\{\theta_{s}\}$, and a starting point  $u^{0}_{i} \in X_{i}$.

 Let us first consider the $k$-th iterate (\ref{eq:4.2}) at $x^{k}=(x^{k}_{i})_{i=1, \ldots, m}$
 within the $s$-th stage. As indicated in (\ref{eq:4.1}), calculation of $x^{k+1}$ then corresponds to the solution of
the optimization problem:
$$
\begin{array}{c}
\displaystyle \min \limits_{z \in X} \rightarrow \left\{
{h(z,e^{s})+\langle p'(x^{k}), z \rangle + (2\alpha )^{-1}
\|z-x^{k}\|^{2}} \right\},
\end{array}
$$
which is equivalent to the $m$ independent problems:
\begin{equation} \label{eq:5.2}
\begin{array}{c}
\displaystyle \min \limits_{z_{i} \in X_{i}} \rightarrow \left\{
{ \varepsilon^{s}_{i} f_{i}(z_{i})+\langle g^{k}_{i}, z_{i} \rangle + (2\alpha )^{-1}
\|z_{i}-x^{k}_{i}\|^{2}} \right\}, \ g^{k}_{i}=\frac{\partial p(x^{k})}{\partial
x_{i}}, \ i=1, \dots, m.
\end{array}
\end{equation}
Note that
\begin{equation} \label{eq:5.3}
g^{k}_{i}=\left\{ {
\begin{array}{ll}
\displaystyle
\tau^{-1} (2x^{k}_{1}-x^{k}_{2}-x^{k}_{m}) \quad & \mbox{if} \ i=1, \\
\tau^{-1} (2x^{k}_{i}-x^{k}_{i+1}-x^{k}_{i-1}) \quad & \mbox{if} \ i=2, \dots, m-1, \\
\tau^{-1} (2x^{k}_{m}-x^{k}_{1}-x^{k}_{m-1}) \quad & \mbox{if} \ i=m, \\
\end{array}
} \right.
\end{equation}
hence each $i$-th agent has the sufficient data for the completely independent
solution of his/her optimization problem (\ref{eq:5.2}).
Afterwards, he/she reports the obtained unique solution $x^{k+1}_{i}$
to the closest neighbours. We suppose that all the units have similar computational capacities
and that the complexity of the private problems (\ref{eq:5.2}) is almost the same.
Therefore, we can suppose that the agents will complete these problems almost simultaneously.

Next, the $i$-th agent should change the current stage in the case where
\begin{equation} \label{eq:5.4}
\|x^{k+1}_{i}- x^{k}_{i} \| \leq \theta_{s}/\sqrt{m},
\end{equation}
cf. (\ref{eq:4.4}). There exist several strategies for
changing the current stage of the method in the whole network that are
dependent of the peculiarities of the system and problem under solution.
If the units and personal problems are very similar to each other, it seems natural to
choose
$$
\varepsilon^{s}_{i}=\sigma_{s}, \ i=1, \dots, m,
$$
for some positive sequence $\{\sigma_{s}\}$. Then the $i$-th agent can in principle
change the current stage of the method for his/her subproblem if the situation (\ref{eq:5.4}) occurs.
In order to adjust the procedure to possible one stage time deviations for different agents
they can send short signals about the current satisfaction of condition (\ref{eq:5.4}).
The simplest protocol is that the $i$-th agent sends such a signal to the closest neighbours in case
(\ref{eq:5.4}), but changes the current stage of the method for his/her subproblem only on receiving
such confirmation signals from both the neighbours.
The more complicated protocol will consists in introducing the same basic
positive sequence $\{\sigma_{s}\}$ and setting
$$
\tilde \varepsilon^{s}_{i}=\sigma_{s}, \ i=1, \dots, m,
$$
and
$$
\varepsilon^{1}_{i}=\tilde \varepsilon^{1}_{i}, \ i=1, \dots, m.
$$
Next, each $i$-th agent also sends the signal to the closest neighbours in case
(\ref{eq:5.4}). On receiving this signal from the $j$-th agent
he/she makes the proper label in his/her list of network agents
and further transmits this signal in the same direction. It is supposed that the
confirmation signals are very short and can be sent at any moment.
After filling out the full list, the $i$-th agent
changes the current stage of the method for his/her subproblem.
In addition, each $i$-th  agent can evaluate the current stage completing
moments for all the agents with taking into account signal delay
time in the network and choose the next number of the sequence $\{\tilde \varepsilon^{q}_{i}\}$
as $\varepsilon^{s+1}_{i}$. For instance, if $\varepsilon^{s}_{i}=\tilde \varepsilon^{t(s)}_{i}$
and his/her $s$-th stage completing time is less essentially than those of most agents,
the $i$-th agent can take $\varepsilon^{s+1}_{i}=\tilde \varepsilon^{t(s)+l(s)}_{i}$, where $l(s) \geq 2$.
Otherwise, $l(s)=1$. This adaptive strategy will equilibrate the stages for
different agents.

The choice of the scaling parameter $\tau > 0$ depends on the desired  value of $\beta \in (0,1)$
in the descent inequality (\ref{eq:2.10}) and the topology of the communication network,
which determines the norm of the matrix $S$ in Lemma \ref{lm:4.1}. In the above setting
we have $\| S\| \leq 4$. Hence, we can provide $\beta=0.5$ if we take e.g.  $\alpha=0.5$ and $\tau=4/3$ or
$\alpha=1$ and $\tau=4$.
Therefore, these conditions give a significant freedom in the choice of the parameters.
In particular, if we choose some other topology of the communication network,
the norm of the matrix $S$ may change, but the choice of the
scaling parameter $\tau$ still provide the suitable
descent property of the forward-backward
splitting method (\ref{eq:4.5}) applied to the penalized problem of form (\ref{eq:3.6}).


\section{Application to the feasibility problem}\label{sc:6}

In the case where $f_{i} \equiv 0$ the original problem
(\ref{eq:1.1})--(\ref{eq:1.2}) reduces to the so-called feasibility problem,
which is to find a point of the set $\tilde X =\bigcap \limits^{m} _{i=1} X_{i}$.
In general, this set may be empty since the set $X_{i}$
may be only known  to the $i$-th agent, then we should find some approximation of
the common point. Following the re-formulation of the original problem given in Section \ref{sc:3}
we obtain the optimization problem
\begin{equation} \label{eq:6.1}
 \min \limits _{x \in  X} \to p(x),
\end{equation}
where the function $p$ is defined in (\ref{eq:3.4}); cf. (\ref{eq:3.4a}).
Clearly, this problem can be solved by the usual gradient projection method,
which corresponds to the splitting method (\ref{eq:2.7}) of Section \ref{sc:2}
 in case $\mu_{1} \equiv 0$ and is written as follows:
\begin{equation} \label{eq:6.2}
\displaystyle x^{k+1}_{i}=\pi_{X_{i}}[x^{k}_{i}-\alpha g^{k}_{i}], \ i=1, \dots, m,
\end{equation}
for $k=0,1,\ldots$, where $g^{k}_{i}$ is defined in (\ref{eq:5.3}). Here and below, $\pi _{V} (v)$ denotes the
projection of $v$ onto $V$. Clearly, the iterates in (\ref{eq:6.2})
are very suitable for the decentralized implementation.
Suppose that the assumptions in (B1) are fulfilled. Then
method (\ref{eq:6.2}) with
$ \alpha \in (0,  \tau/2)$
provides convergence of the sequence
$\{x^{k}\}$, i.e. it has
limit points and all these limit points are solutions of problem (\ref{eq:6.1}).
In fact, due to Lemma \ref{lm:4.1}, $L_{p} \leq 4/\tau$, then the result follows e.g.
from Theorem 1.4 in \cite[Ch.V]{GT89}. The convergence properties of  method (\ref{eq:6.2})
can be also deduced from Proposition \ref{pro:2.5}.
We observe that problem (\ref{eq:6.1}) is solvable under simple sufficient assumptions.
Let each $X_{i}$ be a nonempty, convex, and closed set in
$ \mathbb{R}^{n}$ for $i=1, \dots, m$. Then the set $ X^{*}(p)$ is nonempty if
each $X_{i}$ is a polyhedral set; see \cite{FW56}.
Also, the set $ X^{*}(p)$ is nonempty if at least one of the sets $X_{i}$ is bounded,
then the function $p$ is coercive on $X$.

For many significant applications the above feasibility
problem appears ill-posed, i.e.
its solution does not depend continuously on the input data.
Hence, even small perturbation of the input data may give large deviations from
the solution, which is very essential for the decentralized systems.
In order to overcome these drawbacks, suitable regularization techniques can be applied.
This means that we will again solve problem
(\ref{eq:3.1})--(\ref{eq:3.3})
where
$$
f_{i}(x_{i})=(b\tau )^{-1}\| x_{i}\|^{b}_{d},
 \ d \geq 1,
$$
for $i=1, \dots, m$. Then (B1) implies that the assumptions in (B2) are fulfilled.
Therefore, in accordance with Theorem \ref{thm:4.1}, Method (DPM) will give
the regularized solution of the feasibility
problem.


\section{Computational experiments}\label{sc:7}

In order to check the performance of the proposed methods we carried
out preliminary series of computational experiments.
We chose the basic step for comparison of different methods,
which corresponds to one parallel step of all the $m$ units
and one communication round among the closest neighbours.
Hence, we intend to evaluate the total number of basic steps
of a selected method for obtaining some desired accuracy or a
value of the goal function.
This enables us to utilize a usual PC for the experiments.
We chose two classes of test problems and  implemented all the methods in Delphi with double precision
arithmetic.


\subsection{Feasibility problem tests}\label{sc:7.1}

First we took the affine feasibility problem from Section \ref{sc:6},
which is to find a point of the set
\begin{equation} \label{eq:7.1}
 \tilde X =\left\{v \in \mathbb{R}^{n} \ | \ \tilde Av \leq b \right\},
\end{equation}
where $\tilde A$ is an $m \times n$ matrix, $b\in \mathbb{R}^{m}$.
In other words, we have to solve a system of linear inequalities, which may be inconsistent in general.
We can write $\tilde X =\bigcap \limits^{m} _{i=1} X_{i}$ where
$$
X_{i} =\left\{v \in \mathbb{R}^{n} \ | \ \langle \tilde a_{i}, v\rangle \leq b_{i} \right\},
$$
$\tilde a_{i} $ is the $i$-th row of $\tilde A$, $i=1,\ldots, m$. This means that only the $i$-th unit
of the network knows the vector $\tilde a_{i} $ and number $b_{i} $.
For comparison, we chose the gradient projection method (GPM) given in  (\ref{eq:6.2}),
the classical sequential projection method (SQP) from
 \cite{Agm54} and the alternating direction method (ADM) from
 \cite[Sect. 5.1.2]{BPCPE11}. For (GPM), we fixed the parameters as follows:
 \begin{equation} \label{eq:7.2}
 \alpha =0.4,  \ \tau=1.
 \end{equation}
One iteration of (GPM) or (SQP) will correspond to one basic step,
since one iteration of (SQP) is carried out by one unit, but
all the other units must only wait for the output of the active
unit. Next, one iteration of (ADM) involves parallel calculation of
the current primal points, transmission of these points and
dual points from the previous iteration to a central unit for
the calculation of the average point, which gives $m+1$ basic steps
since one vertex of the communication network may act a central
unit, but it need not be a nearest neighbour to all the units.
Afterwards, the central unit sends the average point to
all the units for updating their dual points. Therefore, we
think that one iteration of (ADM) corresponds to $2m+1$ basic steps.

The standard gap function for problem (\ref{eq:7.1}) is the following:
$$
 \Delta_{s}(v)=\max \limits_{i=1, \dots, m} [\langle \tilde a_{i}, v\rangle - b_{i}]_{+},
$$
where $[\alpha]_{+}=\max\{\alpha,0\}$.
Given a point $x=(x_{i})_{i=1, \ldots, m} \in \mathbb{R}^{N}$, $N=mn$,
we can calculate the average point
\begin{equation} \label{eq:7.4}
 z=(1/m) \sum \limits^{m} _{i=1} x_{i}
 \end{equation}
and then take the standard gap function $\Delta_{s}(z)$
for evaluation of just $x$ in (\ref{eq:6.1}). Besides, we can simply take
the  value
$$
 \Delta_{p}(x)=\sqrt{2p(x)} =\left(\sum \limits^{m-1} _{i=1} \|x_{i}-x_{i+1}\|^{2} + \|x_{m}-x_{1}\|^{2}\right)^{1/2}
$$
as the other gap function for (\ref{eq:6.1}).

\medskip
\noindent {\em Example 1 (Consistent case).}
The elements of the matrix $\tilde A$ were defined by
$$
\tilde a_{ij}= \left\{
\begin{array}{rl}
\displaystyle
-0.2 ij \quad & \mbox{for} \ j=1, \dots, n/2, \\
0.2 ij \quad & \mbox{for} \ j=n/2+1, \dots, n,
\end{array}
 \right. \quad \mbox{if} \ i \ \mbox{is odd},
$$
and
$$
\tilde a_{ij}= \left\{
\begin{array}{rl}
\displaystyle
0.2 (i-1)(n+1-j) \quad & \mbox{for} \ j=1, \dots, n/2, \\
-0.2 (i-1)(n+1-j) \quad & \mbox{for} \ j=n/2+1, \dots, n,
\end{array}
 \right. \quad \mbox{if} \ i \ \mbox{is even},
$$
elements of the vector $b$ were defined by
$$
b_{i}=\sum_{j=1}^{n} \tilde a_{ij} \ \mbox{ for } \ i=1, \dots, m,
$$
where $m$ and $n$ were chosen to be even and $m > n$.
It follows that the system is solvable and involves in fact only two
different inequalities. Due to the distributed treatment of the problem,
this simple example can be chosen for calculations.
We took the same starting point $x^{0} = (5, \dots,5)^{\top} $
and the accuracy $\delta=0.0001$ with respect to $\Delta_{p}(x)$ for all the methods.

Table \ref{tbl:7.1} describes the results of application of (GPM) to this problem,
where (kt) denotes the number of the basic iterations. The third column gives the total
number of the basic iterations for attaining the accuracy $\delta$ with respect to $\Delta_{p}(x)$.
Columns 4--6 show the attained values of $\Delta_{p}(x)$ and $\Delta_{s}(z)$ for kt=10,20,30.
Note that the average point $z$ in (\ref{eq:7.4}) was calculated in a separate block only for
derivation of the current value of $\Delta_{s}(z)$, and these values were not used for
the method itself.
\begin{table}
\caption{Example 1 for (GPM)} \label{tbl:7.1}
\begin{center}
\begin{tabular}{|rr|r|rr|rr|rr|}
\hline
       &     &     &       &      kt=10       &       &  kt=20          &     &  kt=30   \\
\hline
   $m$ & $n$ &  kt & $\Delta_{p}$ & $\Delta_{s}$  & $\Delta_{p}$ & $\Delta_{s}$  & $\Delta_{p}$ & $\Delta_{s}$  \\

\hline
$20$ & $10$  & 32   & 0.81 & 5.97 & 0.01  & 0.09 &  0.0002 & 0.0014  \\
\hline
$50$ & $10$  & 33   & 1.28 & 15.39 & 0.02 & 0.23  &   0.0003 & 0.0035  \\
\hline
$100$ & $10$ & 34   & 1.81 & 31.11 & 0.03 & 0.47  &  0.0004 & 0.007   \\
\hline
$100$ & $20$  & 32 & 1.77 & 82.37 & 0.02 & 0.81  &   0.0002 & 0.0018  \\
\hline
$100$ & $50$  & 31 & 2.21 & 334.42 & 0.02 & 2.9  &  0.0001 & 0.0214  \\
\hline
\end{tabular}
\end{center}
\end{table}
Table \ref{tbl:7.2} describes the results of application of (SQP),
where the third column gives the total
number of the basic iterations for attaining the accuracy $\delta$ with respect to $\Delta_{p}(x)$.
Columns 4--6 show the attained values of $\Delta_{s}(x)$ for kt=10,20,30.
Note that they were calculated at the current iterates
rather than at some average points.
\begin{table}
\caption{Example 1 for (SQP)} \label{tbl:7.2}
\begin{center}
\begin{tabular}{|rr|r|r|r|r|}
\hline
       &     &     &    kt=10       &  kt=20          &   kt=30   \\
\hline
   $m$ & $n$ &  kt &  $\Delta_{s}$  &  $\Delta_{s}$  &  $\Delta_{s}$  \\

\hline
$20$ & $10$  & 40    & 3.26 &  0.01 &  0.0001   \\
\hline
$50$ & $10$  & 100   & 8.41 &  0.03  &   0.0001  \\
\hline
$100$ & $10$ & 200   & 17 &  0.06  &  0.0002    \\
\hline
$100$ & $20$  & 200  & 37.58 &  0.07  &   0.0003   \\
\hline
$100$ & $50$  & 200  & 159.14 &  0.21  &  0.0003  \\
\hline
\end{tabular}
\end{center}
\end{table}
Table \ref{tbl:7.3} describes the results of application of (ADM),
where (kl) denotes the number of its iterations. Columns 3--4 show the total
number of both the basic and its iterations for attaining the accuracy $\delta$ with respect to $\Delta_{p}(x)$.
Columns 5--7 show  the
number of these iterations for attaining a solution accuracy $\delta$ with respect to $\Delta_{s}(z)$.
\begin{table}
\caption{Example 1 for (ADM) } \label{tbl:7.3}
\begin{center}
\begin{tabular}{|rr|rr|rrr|}
\hline
   $m$ & $n$ &  kt  & kl  &  kt & kl  &  $\Delta_{s}$  \\

\hline
$20$ & $10$  & 287  & 7   &  144  & 3 &  0   \\
\hline
$50$ & $10$  & 707  & 7   &  354  & 3 &  0  \\
\hline
$100$ & $10$ & 1407  & 7   &  704  & 3 &  0     \\
\hline
$100$ & $20$  & 1407  & 7   &  704  & 3 &  0    \\
\hline
$100$ & $50$  & 1407  & 7   &  704  & 3 &  0  \\
\hline
\end{tabular}
\end{center}
\end{table}

\medskip
\noindent {\em Example 2 (Inconsistent case).}
The elements of the matrix $\tilde A$ were defined by
$$
\tilde a_{ij}=  2 \sin(i/j) \cos(ij), \quad  j=1, \dots, n,  \quad \mbox{if} \ i =1, \dots, m, \ i \neq n,
$$
and
$$
\tilde a_{nj} =-\sum_{i=1}^{n-1} \tilde a_{ij} \quad  j=1, \dots, n;
$$
elements of the vector $b$ were defined by
$$
b_{i}=
\left\{
\begin{array}{rl}
\displaystyle
\sum_{j=1}^{n} \tilde a_{ij} -5, \quad & i=1, \dots, n, \\
\displaystyle
\sum_{j=1}^{n} \tilde a_{ij} +5, \quad & \ i=n+1, \dots, m,
\end{array}
\right.
$$
where $m$ and $n$ were chosen to be even and $m > n$.
It follows that the system in (\ref{eq:7.1}) is inconsistent, hence we have to introduce
the other basic gap function:
$$
 \Delta_{d}(x)=\|x-\pi_{X}[x-\alpha p'(x)]\|,
$$
where $\alpha$ and $\tau$ are given in (\ref{eq:7.2}).
We took the same starting point $x^{0} = (5, \dots,5)^{\top} $ for all the methods.
Table \ref{tbl:7.4} describes the results of application of (GPM) to this problem.
Column 3 shows the numbers of the basic iterations (kt) for attaining
the accuracy $\delta=0.1$ with respect to $\Delta_{d}(x)$.
Columns 4--5 show the values of the basic iterations
for attaining the accuracy $\delta=0.01$ and the related values of
$\Delta_{p}(x)$. Column 6 shows  the values of $\Delta_{s}(z)$ and
the corresponding numbers of the basic iterations.
\begin{table}
\caption{Example 2 for (GPM)} \label{tbl:7.4}
\begin{center}
\begin{tabular}{|rr|r|rr|l|}
\hline
       &     &  $\delta=0.1$   &              &   $\delta=0.01$       &         \\
\hline
   $m$ & $n$ &  kt             & kt           & $\Delta_{p}$          & kt / $\Delta_{s}$  \\

\hline
$20$ & $10$  & 108   & 597 & 6.46 & 580/12.25    \\
\hline
$50$ & $10$  & 93   & 897 & 6.31 & 880/12.07   \\
\hline
$100$ & $10$ & 125   & 836 & 6.34 & 820/8.98   \\
\hline
$100$ & $20$  & 220   & 2176 & 4.14 & 2160/10.03  \\
\hline
$100$ & $50$  & 280   & 5038 & 3.06 & 5020/192.67   \\
\hline
\end{tabular}
\end{center}
\end{table}
Table \ref{tbl:7.5} describes the results of application of (SQP) and  (ADM) to this problem.
Both the methods do not converge. We give only one experiment for (SQP), where $\Delta_{d}$
indicated the average distance after one cycle of $m$ iterations:
$$
 \Delta_{d}(x)=(1/m)\left(\sum_{k=1}^{m-1} \|x^{k+1}-x^{k}\|+\|x^{m}-x^{1}\|\right).
$$
The results of application of (ADM) are given for all the variants,
where $\Delta_{d}$
indicated the distance between two average primal points:
$$
 \Delta_{d}(z)= \|z^{k}-z^{k-1}\|.
$$
\begin{table}
\caption{Example 2 for (SQP) and  (ADM)} \label{tbl:7.5}
\begin{center}
\begin{tabular}{|rr|rr|rrr|}
\hline
       &     &         &       &              &        &    (SQP)     \\
\hline
   $m$ & $n$ &  kt    & kl    & $\Delta_{d}$   & $\Delta_{p}$    & $\Delta_{s}$  \\
\hline
$20$ & $10$  & 1000   & -     & 0.37 & 2.78 & 12.09    \\
\hline
       &     &         &       &              &        &    (ADM)     \\
\hline
   $m$ & $n$ &  kt    & kl    & $\Delta_{d}$   & $\Delta_{p}$    & $\Delta_{s}$  \\
\hline
$20$ & $10$  & 1005   & 24     & 0.56 & 8.82 & 9.74     \\
\hline
$50$ & $10$  & 960   & 10     & 1.29 & 11.14 & 12.77   \\
\hline
$100$ & $10$ & 1910   & 9     & 1.3 & 13.09 & 9.8  \\
\hline
$100$ & $20$  & 4925   & 24     & 0.95 & 10.72 & 20.67  \\
\hline
$100$ & $50$  & 9950   & 49     & 1.07 & 14.35 & 25.44 \\
\hline
\end{tabular}
\end{center}
\end{table}
We can conclude that both (SQP) and  (ADM) show rather rapid convergence for
consistent systems, but they are not adjusted for the decentralized implementation
since the iteration points may be far from each other unlike (GPM).
Besides, (GPM) shows rather stable convergence for
inconsistent systems, this is not the case for (SQP) and  (ADM).


\subsection{Fermat-Weber problem tests}\label{sc:7.2}

We also took the well-known Fermat-Weber problem in the following simple format:
$$
 \min \limits _{v \in \mathbb{R}^{n}} \to \varphi(v)=\sum \limits^{m} _{i=1} \|v-\tilde a_{i}\|,
$$
where $\tilde a_{i} $, $i=1,\ldots, m$ are some given points (anchors). Clearly, this is a
particular case of problem (\ref{eq:1.1})--(\ref{eq:1.2}) with the coercive, convex,
and non-smooth cost function $\varphi$ over the whole space $\mathbb{R}^{n}$. It can be rewritten in the format
(\ref{eq:3.1})--(\ref{eq:3.3}) as follows:
$$
 \min \limits _{x \in Y} \to \sum \limits^{m} _{i=1} \|x_{i}-\tilde a_{i}\|,
$$
where
$$
Y=\left\{x \in \mathbb{R}^{N} \ | \ x_{i}=x_{i+1}, \ i=1, \dots, m-1, x_{m}=x_{1}  \right\},
$$
i.e. $X_{i}=\mathbb{R}^{n}$, $f_{i}(x_{i})=\|x_{i}-\tilde a_{i}\|$ for $i=1,\ldots, m$.
Then its solution set  $D^{*}(f)$ is nonempty
and  bounded. In the multi-agent setting,  the $i$-th unit
of the network knows only the vector $\tilde a_{i} $.

For comparison, we chose the two-level decomposable penalty method  (DPM) from
Section \ref{sc:4} and the primal-dual method (PDM) with the proper adjustment to
the same multi-agent network setting; see \cite{CP16,AH16,LLZ20}.

The $k$-th iteration of (PDM), $k=0,1,2,\ldots$, is described as follows.
Given a pair $(x^{k},w^{k})$,
$x^{k}=(x^{k}_{i})_{i=1, \ldots, m} \in \mathbb{R}^{N}$,
$w^{k}=(w^{k}_{i})_{i=1, \ldots, m} \in \mathbb{R}^{N}$, the units independently find
points $x^{k+1}_{i}$ as solutions of the optimization problems:
$$
\begin{array}{c}
\displaystyle \min \limits_{z_{i} } \rightarrow \left\{
{ \|z_{i}-\tilde a_{i}\| +\langle \tilde g^{k}_{i}, z_{i} \rangle + (2\alpha )^{-1}
\|z_{i}-x^{k}_{i}\|^{2}} \right\},
\end{array}
$$
where
$$
\tilde g^{k}_{i}=\left\{ {
\begin{array}{ll}
\displaystyle
w^{k}_{1}-w^{k}_{m} \quad & \mbox{if} \ i=1, \\
w^{k}_{i}-w^{k}_{i-1} \quad & \mbox{if} \ i=2, \dots,m.
\end{array}
} \right.
$$
Then they calculate the points $\tilde x^{k+1}_{i}=2x^{k+1}_{i}-x^{k}_{i}$
for $i=1,\ldots, m$, report these points
to the neighbours, and find dual
points $w^{k+1}_{i}$ as solutions of the optimization problems:
$$
\begin{array}{c}
\displaystyle \min \limits_{u_{i} } \rightarrow \left\{
{(2\beta )^{-1} \|u_{i}-w^{k}_{i}\|^{2} -\langle \tilde q^{k}_{i}, u_{i} \rangle  } \right\},
\end{array}
$$
where
$$
\tilde q^{k}_{i}=\left\{ {
\begin{array}{ll}
\displaystyle
\tilde x^{k+1}_{i}-\tilde x^{k+1}_{i+1} \quad & \mbox{if} \ i=1, \dots,m-1, \\
\tilde x^{k+1}_{m}-\tilde x^{k+1}_{1} \quad & \mbox{if} \ i=m.
\end{array}
} \right.
$$
Then they  report these points to the neighbours.

Note that, unlike (DPM), each $i$-th agent calculates both primal and dual points and
twice reports the obtained points
to the neighbours.  Hence, one iteration of (PDM) corresponds to two basic steps.

For (DPM), we fixed the parameters $\alpha$ and $\tau$ as in (\ref{eq:7.2}), besides,
 we used the rule
 $$
\theta_{s+1}=q_{1} \theta_{s}, \ \sigma_{s+1}=q_{2} \sigma_{s}, \ \varepsilon^{s}_{i}=\sigma_{s}, \ i=1, \dots, m, \
0<q_{1} <q_{2} <1,
$$
where $\theta_{0}=0.5$, $\sigma_{0}=1$. For (PDM), we fixed its parameters as follows:
$$
 \alpha =0.5,  \ \beta=0.25.
$$
We took the same starting point $x^{0} = (5, \dots,5)^{\top} $ for both the methods.
Since the methods are essentially different we took the value of the cost function $\varphi(z)$
calculated at the average primal point $z$ from (\ref{eq:7.4}) after the same number of
basic steps. This average point was calculated in a separate block
 and this value was not used in the methods themself.

\medskip
\noindent {\em Example 3 (Exact case).}
The elements of the vectors $\tilde a_{i}$ were defined by
$$
\tilde a_{ij}=  5 \sin(i/j) \cos(ij), \quad  j=1, \dots, n, \ i =1, \dots, m.
$$
The parameters $q_{1}$ and $q_{2}$ in (DPM) were chosen as follows:
$$
q_{1}=0.1, \ q_{2}=0.6.
$$
Table \ref{tbl:7.6}  shows the values of $\varphi(z)$ for various
numbers of the basic iterations of (DPM) and  (PDM).
\begin{table}
\caption{Example 3 for (DPM) and  (PDM)} \label{tbl:7.6}
\begin{center}
\begin{tabular}{|r|rr|rrrr|}
\hline
 (DPM) &  $m$ & $n$ &  kt=0    &  kt=60     &  kt=100    &  kt=200   \\
\hline
 & $20$ & $10$  & 360.85   & 155.82    & 152.6 & 152.36     \\
\hline
 & $50$ & $10$  & 875.72   & 388.64     & 382.82 & 382.28   \\
\hline
 & $100$ & $10$ & 1747.73   & 771.74     & 760.17 & 759.42  \\
\hline
 & $100$ & $20$  & 2495.44   & 1197.44     & 1100.81 & 1095.09   \\
\hline
&  $100$ & $50$  & 3951.23   & 2373.52     & 1902.42 & 1764.77 \\
\hline
 (PDM) &  $m$ & $n$ &  kt=0    &  kt=60     &  kt=100    &  kt=200   \\
\hline
 & $20$ & $10$  & 360.85   & 181.08    & 155.14 & 152.34     \\
\hline
 & $50$ & $10$  & 875.72   & 443.04     & 388.33 & 382.25   \\
\hline
 & $100$ & $10$ & 1747.73   & 880.19     & 771.53 & 759.41  \\
\hline
 & $100$ & $20$  & 2495.44   & 1492.05     & 1193.02 & 1096.12   \\
\hline
&  $100$ & $50$  & 3951.23   & 2871.7     & 2343.01 & 1816.31 \\
\hline
\end{tabular}
\end{center}
\end{table}
We also noted that  replacing rule (\ref{eq:4.4}) with  (\ref{eq:5.4}) in (DPM) did not affect the convergence.

\medskip
\noindent {\em Example 4 (Perturbed case).}
The vectors $\tilde a_{i} $, $i=1,\ldots, m$ were defined as in Example 3.
The main difference was in inserting perturbations in transmitted data.
That is, all the neighbours of the $i$-th unit in (DPM) received the perturbed value
$x_{ij}+0.5 \sin(i) \sin(j)$ instead of $x_{ij}$. The same data perturbations were inserted
in (PDM) for $\tilde x_{ij}$ and $w_{ij}$.
The parameters $q_{1}$ and $q_{2}$ in (DPM) were chosen as follows:
$$
q_{1}=0.2, \ q_{2}=0.5.
$$
Table \ref{tbl:7.7}  shows the values of $\varphi(z)$ for various
numbers of the basic iterations of (DPM) and  (PDM) in this case.
\begin{table}
\caption{Example 4 for (DPM) and  (PDM)} \label{tbl:7.7}
\begin{center}
\begin{tabular}{|r|rr|rrrr|}
\hline
 (DPM) &  $m$ & $n$ &  kt=0    &  kt=60     &  kt=100    &  kt=200   \\
\hline
 & $20$ & $10$  & 360.85   & 156.1    & 153 & 152.59    \\
\hline
 & $50$ & $10$  & 875.72   & 388.64     & 383.12 & 382.36   \\
\hline
 & $100$ & $10$ & 1747.73   & 771.73     & 760.44 & 759.5  \\
\hline
 & $100$ & $20$  & 2495.44   & 1197.4     & 1100.93 & 1095.36   \\
\hline
&  $100$ & $50$  & 3951.23   & 2373.65     & 1902.53 & 1765.63 \\
\hline
 (PDM) &  $m$ & $n$ &  kt=0    &  kt=60     &  kt=100    &  kt=200   \\
\hline
 & $20$ & $10$  & 360.85   & 181.19    & 155.18 & 152.35     \\
\hline
 & $50$ & $10$  & 875.72   & 443     & 388.32 & 382.26   \\
\hline
 & $100$ & $10$ & 1747.73   & 880.14     & 771.52 & 759.41  \\
\hline
 & $100$ & $20$  & 2495.44   & 1492.06     & 1193.03 & 1096.12   \\
\hline
&  $100$ & $50$  & 3951.23   & 2871.7     & 2343 & 1816.3 \\
\hline
\end{tabular}
\end{center}
\end{table}

We can conclude that small data perturbations did not affect significantly
the convergence of (DPM) and  (PDM). In the first period of work, (DPM) appeared more rapid,
i.e. it is better for obtaining some good approximation of the solution.
Afterwards, its convergence appeared somewhat slower in comparison with (PDM), but
for large dimensionality (DPM) still had some preference.


\section{Conclusions}\label{sc:8}

We described a new decentralized penalty method for
convex constrained optimization problems
in a decentralized multi-agent network setting.
Its convergence was established under rather weak assumptions,
even if the constraints were inconsistent. The computational
experiments confirmed rather satisfactory convergence.
Nevertheless, there are several directions for further investigations.
In particular, they involve applications to different
computational network topologies and to different classes of multi-agent
optimization problems. Also, proper choice of the parameters with respect to
special problems needs additional substantiation.


\section*{Acknowledgement}

This work was supported  by grant No. 331833 from Academy of Finland
and by the RFBR grant, project No. 19-01-00431.


\end{document}